\documentclass[11pt]{article}
\usepackage{enumerate}
\usepackage{amssymb}
\usepackage{amsmath,amsfonts,amsthm}
\usepackage[colorlinks=true]{hyperref}
\usepackage[english]{inputenc}

\numberwithin{equation}{section}

\newtheorem{theorem}{Theorem}[section]

\newtheorem{prop}{Proposition}[section]
\newtheorem{lem}{Lemma}[section]
\newtheorem{rmk}{Remark}[section]

\def\nd{\noindent}

\def\w{W_0^{1,\Phi}(\Omega )}

\def\nd{\noindent}

\def\Re{\mathbb{R}}

\begin{document}

\title{ Sign changing solutions for quasilinear superlinear\\
\nd  elliptic problems}
\vskip.5cm

\author{M. L. M. Carvalho ~-~F. J.  S. A.  Corr\^ea ~-~ Jose V. A. Goncalves~-~E. D.  Silva}

\vskip.5cm

\date{}

\pretolerance10000

\maketitle

\nd  {\bf Abstract} Results on existence and multiplicity of solutions for a nonlinear elliptic problem driven  by the $\Phi$-Laplace operator are established. We employ minimization arguments on suitable Nehari manifolds to build  a negative and a positive ground state solutions. In order to find a nodal solution we employ additionally the well known Deformation Lemma and Topological Degree Theory.
\vskip1cm

 \nd {\bf Keywords}~ {Variational methods~$\cdot$~Quasilinear Elliptic Problems~$\cdot$~ Nehari manifold method}
\vskip.2cm

\nd  {\bf Mathematics Subject Classification (1991)} {Primary~ 35J20 $\cdot$  35J25 $\cdot$ 35J60}
\vskip.7cm

\section{\Large \bf  Introdution}

\nd In this work we consider the quasilinear elliptic problem

\begin{equation}\label{pi}
\left\{
\begin{array}{l}
- {div}\left( \phi(|\nabla u|) \nabla u \right) = f(x,u) \
\mbox{in} \ \ \Omega, \ \ \\
  u = 0 \ \mbox{on} \ \
\partial\Omega,
\end{array}
 \right.
\end{equation}

\nd where $\Omega\subset\mathbb{R}^{N}$ is a  bounded,  smooth domain, $f: \Omega \times  \mathbb{R} \rightarrow \mathbb{R}$ is a Caratheodory function which is  $C^{1}$ in the second variable.  For the function $\phi$ we assume that $\phi: (0,\infty)\rightarrow (0,\infty)$ is of class $C^{2}$   and satisfies the following conditions:
\vskip.1cm

\begin{itemize}
  \item[$(\phi_1)$] $\displaystyle \lim_{t \rightarrow 0} t \phi(t)= 0, \displaystyle \lim_{t \rightarrow \infty} t \phi(t)= \infty$;
  \item[$(\phi_2)$] $t\mapsto t\phi(t)$ is strictly increasing.
\end{itemize}

\nd We point out that the function  $\phi(t) = t^{p-2}~\mbox{for}~t > 0$, with $1< p < \infty$, satisfies  $(\phi_{1})- (\phi_{2})$ and in this case the operator in  problem (\ref{pi}) is named $p$-Laplacian and (\ref{pi}) reads as
\begin{equation*}\label{ar}
\displaystyle  -\Delta_p u = f(x,u)~ \mbox{in}~
\Omega,~~ u = 0~  \mbox{on}~ \partial \Omega.
\end{equation*}
\vskip.1cm

\nd  In a similar  way, the function  $\phi(t) = t^{p -2} + t^{q -2}$ with $ 1 < q < p < \infty$ satisfies the conditions $(\phi_{1})-(\phi_{2})$. In this case the operator in  problem (\ref{pi}) is named $(p,q)$-Laplacian  and problem \eqref{pi} becomes
\begin{equation*}\label{arr}
\displaystyle  -\Delta_p u - \Delta_q u = f(x,u)~ \mbox{in}~
\Omega,~~ u = 0~  \mbox{on}~ \partial \Omega.
\end{equation*}

\nd {---------------------------}

\nd F. J. S. A.  Corr\^ea and E. D.  Silva   were supported in part  by CNPQ/Brazil.

\nd {---------------------------}


\nd In the model case $f(s)=|s|^{p-2}s$~  it is well-known  that the Ambrosetti-Rabinowitz conditon  (see \cite{AR}),  $(AR)$ for short,  namely
$$
\mbox{\it there exist}~  \theta > 2 , R > 0~  \mbox{\it such that}
$$
\begin{equation} \tag{$AR$}\label{AR}
0  <\theta F(x,t) \leq t f(x,t) , x \in \Omega, |t| \geq R,
\end{equation}

\nd with $F(x,t)=\displaystyle\int_{0}^{t}f(x,s)ds$, plays  a crucial role while   addressing  compactness requirement in variational methods.  However, there exist  lots  of functions for which \eqref{AR} is not satisfied. For instance, $f(t) = t ~ \log (1 + |t|)~\mbox{for}~ t \in \mathbb{R}$ does not satisfy \eqref{AR}. It is important to emphasize that the main role of $(AR)$  is to ensure the well known  $(PS)$-condition required by  minimax arguments.
\vskip.1cm

 \nd  We refer the reader to the  reasearch papers  \cite{chung, LamLu,  Liu, LiWangZheng, MiyagakiSouto,MugnaiPapageorgiou, fang} and  references therein, where problems involving  the $p$-Laplacian, sometimes  with $p=2$,  and  the $(p,q)$-Laplacian operator have been addressed.
\vskip.1cm

\nd There is a rich literature on problems of the form (\ref{pi}) with functions $\phi$ even more general than the ones mentioned above. In such more general settings the operator in (\ref{pi}) is named $\Phi$-Laplacian and  is written  as
$$
  {div} ( \phi(|\nabla u|) \nabla u) := \Delta_{\Phi} u,
$$
\nd with
$$
\Phi(t) = \int_{0}^{t} s \phi(s) ds, t \in \mathbb{R},
$$
where  $0<s \rightarrow s\phi (s)$ has been extended to the whole $\mathbb{R}$ as an odd function and so $\Phi$ is an even function.
In  \cite{hui-2},  Cl\'ement, Garc\'ia-Huidobro, Man\'asevich \&  Schmitt   showed  results for problems of the form
$$
- \Delta_{\Phi} u = f(x,u)~\mbox{in} \ \ \Omega,~~ u = 0 \ \mbox{on}~\partial\Omega,
$$
\nd  where the nonlinear term $f$ satisfies (AR).  In the recent paper \cite{JVMLED},  results on existence and multiplicity of solutions were proven.

\nd We refer the reader to Radulescu \cite{radu-1} and its references where operators even more general than the $\Phi$-Laplacian are treated, and motivation from the physical sciences  are discussed.
\vskip.2cm

\nd Due to the nature of the operator $\Delta_{\Phi} $  we shall work in the framework of Orlicz-Sobolev spaces $\w$.  Some basic facts and references on these spaces are given in section \ref{O-S-1}.
\vskip.3cm

\nd We shall assume  the following condition on $\phi$:
\vskip.3cm

\begin{description}
\item[$(\phi_3)$]~~~~~~~~ $\displaystyle  -1 <  \ell - 2  := \inf_{t > 0} \dfrac{(t \phi(t))^{\prime \prime }t}{(t \phi(t))^{\prime}} \leq \sup_{t > 0} \dfrac{(t \phi(t))^{\prime \prime}t}{(t \phi(t))^{\prime}} : = m -2 < N - 2.$
\end{description}
\vskip.3cm

\begin{rmk}\label{Reflexivity}
 \nd   It can be shown  that $(\phi_3)$ implies the  (less restrictive) condition:
    \begin{description}
 \item[$(\phi_3)^{\prime}$]~~~~~~~~ $\displaystyle 1<\ell :=\inf_{t>0}\frac{t^{2} \phi(t)}{\Phi(t)}\leq \sup_{t>0}\frac{t^{2} \phi(t)}{\Phi(t)}=:m <N$,
\end{description}
\nd Moreover under conditions $(\phi_{1}), (\phi_{2}), (\phi_{3})^{'}$ the space $W^{1,\Phi}_{0}(\Omega)$ is a reflexive Banach space, (see Remark $\ref{Delta_2}$  at  Section \ref{O-S-1} ahead).
\end{rmk}
\vskip.3cm

\nd We denote by  $\lambda_{1} > 0$ the first eigenvalue for the operator $-\Delta_{\Phi}$.  Recall that it satisfies the Poincar\'e inequality, (see e.g.  \cite{hui-2}, \cite{gossez-Czech}),
\begin{equation}\label{lambda1}
\lambda_{1} \int_\Omega \Phi(u) dx \leq \int_\Omega \Phi(|\nabla u|)dx,~u \in  \w.
\end{equation}

\nd The following conditions will be imposed  on the nonlinear term $f$:

\begin{itemize}\label{f0}
  \item [$(f_0)$] there are a  function  $\psi: [0, \infty) \rightarrow [0, \infty)$ and a constant $C > 0$ such that
  $$|f(x, t)| \leq C \left[ 1 + \psi(t)\right],  \quad  t \in \mathbb{R},~ x \in \Omega,$$
\end{itemize}
where $\Psi(t) =\displaystyle \int_{0}^{t} \psi(s) ds$ is an $N$-function satisfying  $\Psi<<\Phi_*$,  (definitions and properties  in Section \ref{O-S-1}),  and
\begin{description}
\item[($\psi_1$)]~~~~~$\displaystyle 1<\ell\leq m<\ell_\Psi:=\inf_{t>0}\frac{t\psi(t)}{\Psi(t)}\leq
\sup_{t>0}\frac{t\psi(t)}{\Psi(t)}=:m_\Psi<\ell^*:=\frac{\ell N}{N-\ell}$;
\end{description}

\begin{itemize}\label{f1}
  \item [$(f_1)$] the function
 $$t\mapsto  \frac{f(x, t)}{|t|^{m-2}t}$$
 is increasing on   $\mathbb{R}\backslash \{ 0\}$;
  \end{itemize}

\begin{itemize}\label{f2}
  \item [$(f_2)$] the limit
  $$\lim_{t \rightarrow 0} \frac{f(x, t)}{t \phi(t)} < \lambda_{1}$$
  holds uniformly in $x \in \Omega$;
\end{itemize}

\begin{itemize}\label{f3}
\item [$(f_3)$] the limit
$$\lim_{|t| \rightarrow \infty} \frac{f(x, t)}{|t|^{m-2}t} = + \infty$$
holds uniformly for  $x \in \Omega$.
\end{itemize}

 \nd The key technique in the proof of the main result of the present paper will be the use of the  Nehari manifold method, (see the pioneering work \cite{Nehari}),  this time for the energy functional J associated with problem (\ref{pi}) which due to the nature of the operator $\Delta_{\Phi}$ is defined in an Orlicz-Sobolev space. Here one of the main difficulties is to ensure that a minimizing sequence for J over the Nehari  manifold converges to a critical point. Another difficulty  is to show that  the Nehari manifold is  $C^{1}$. To overcome these difficulties we prove that the map $u \mapsto \langle J^{\prime}(u), u\rangle$ is of class $C^{1}$. Condition $(\phi_{3})$ is crucial to achieve these steps.  We will also make use of some Topological Degree arguments.
\vskip.2cm

\nd The main result in this work is:

\begin{theorem}\label{th2}
Assume  $(\phi_{1}), (\phi_{2}), (\phi_{3})$, $(f_{0}) - (f_{3})$. Then  problem \eqref{pi} admits at least two ground state solutions $u_{1}, u_{2} \in W^{1, \Phi}_{0}(\Omega)$ satisfying $u_{1} < 0$ and $u_{2} > 0$ in $\Omega$. In addition,  problem \eqref{pi} admits a further solution say $u_{3}$ which changes sign in $\Omega$.
\end{theorem}

\nd To our best knowledge, there is no result on  existence of sign changing solutions for  problems involving the $\Phi$-Laplacian operator.
\vskip.1cm

\nd Existence of positive and negative solutions have long been searched for problems involving both  Laplacian and  p-Laplacian  equations. We would like to mention  the works  \cite{JVMLED}, \cite{hui-2}, \cite{chung},  \cite{Fuk_1}, \cite{Fuk_2},  \cite{iacopetti}, \cite{fang},  \cite{Wang} and references therein which are more related to our present interest in this paper.  In those works the authors have used  truncation techniques and maximum principles.

\vskip.1cm

\nd  More recently,  sign changing solutions have been considered. We refer the reader to Szulkin and  Weth \cite{SW,ST}, where the authors also  addressed  existence of positive solutions for the Dirichlet problem for  Laplace and p-Laplace equations  using the Nehari manifold method.

\nd We also refer the reader to \cite{TT}, \cite{TTZ}, \cite{Drabek} and \cite{Zou}. In these works the authors considered
semilinear/superlinear problems driven by  the Laplace and the p-Laplace operators. In \cite{TTZ} the authors considered quasilinear  problems for the  p-Laplacian operator obtaining existence of positive solutions. The results in the present paper for the case of the  $\Phi$-Laplace operator  complement/extend ones in the   above-mentioned papers.

\nd We also point out that quasilinear elliptic problems have been considered under several assumptions on the nonlinear term f. In this regard we refer the reader to \cite{JVMLED,chung,gossez-Czech,LamLu,fang, cerami}. In \cite{chung}, \cite{LamLu} the authors considered  monotonocity conditions on the nonlinear term f, proving existence of positive solutions.  In \cite{JVMLED}, \cite{LamLu} the authors studied  quasilinear elliptic  problems without the Ambrosetti-Rabinowitz condition at infinity obtaining existence of multiple  solutions.

\begin{rmk}
 An example of a problem to which our theorem {\rm \ref{th2}} applies is
$$
-\Delta_{\Phi} u=pu^{p-1}\log(1+u)+\frac{u^p}{u+1}~\mbox{in}~\Omega,\\
\\
u\geq0~\mbox{in}~\Omega,~~u=0~~\mbox{on}~\partial \Omega,
$$
\nd where
$$
\Phi(t)=|t|^\gamma\log(1+|t|)~~ \mbox{with}~~ 1<\frac{-1+\sqrt{1+4N}}{2}<\gamma<N-1,
$$
$$
F(t)=t^p\log(1+|t|)~~\mbox{and}~~ F^{'}(t)= p t^{p-1}\log(1+t)+\frac{t^p}{t+1} :=   f(t),~t > 0 ,
$$
\nd with   $\gamma = \ell$,  $\gamma+1 = m$,  $\sigma > \frac{N}{\ell}$~ and~  $m<p<\frac{\ell\sigma}{\sigma-1}$. We notice that    by an easy computation,
$$
 t\phi(t) = \Phi^{\prime}(t) = \gamma t^{\gamma-1}\log(1+t)+\frac{t^\gamma}{1+t},~t > 0
$$
\nd \nd and the conditions $(\phi_1), (\phi_2)$,  $(\phi_3)$ are satisfied.   In addition, the function $f(t)$  satisfies $(f_0)-(f_3)$.
\vskip.2cm

\nd The operator $\Delta_\Phi$ in the present example appears in Plasticity, see e.g. Fukagai  and Narukawa \mbox{\rm \cite{Fuk_2}}.
\end{rmk}

\nd This paper is organized as follows: in section 2 we recall some basic properties of Orlicz-Sobolev spaces.  Section 3 is devoted to  auxiliary results on functionals defined on Orlicz-Sobolev spaces and  related Nehari manifolds. In Section 4 we give the proof of the main result of Section 3, namely Theorem \ref{th1} which ensures existence of a ground state solution of problem (\ref{pi}). Finally Section 5 is devoted to the proof of Theorem \ref{th2}.

\section{\bf \Large Basics on Orlicz-Sobolev spaces}\label{O-S-1}

The reader is  referred to  \cite{A,Rao1} regarding Orlicz-Sobolev spaces.  The usual norm on $L_{\Phi}(\Omega)$    (Luxemburg norm) is ,
$$
\|u\|_\Phi=\inf\left\{\lambda>0~|~\int_\Omega \Phi\left(\frac{u(x)}{\lambda}\right) dx \leq 1\right\}
$$
and  the  Orlicz-Sobolev norm of $ W^{1, \Phi}(\Omega)$ is
\[
  \displaystyle \|u\|_{1,\Phi}=\|u\|_\Phi+\sum_{i=1}^N\left\|\frac{\partial u}{\partial x_i}\right\|_\Phi.
\]
Recall that
$$
\widetilde{\Phi}(t) = \displaystyle \max_{s \geq 0} \{ts - \Phi(s) \},~ t \geq 0.
$$

\begin{rmk}\label{Delta_2}
 We notice that, $(\phi_1), (\phi_2),(\phi_3)^{\prime}$  imply that  $\Phi$ and $\widetilde{\Phi}$  are  N-functions  satisfying  the $\Delta_2$-condition.  In addition,   $L_{\Phi}(\Omega)$  and $W^{1,\Phi}(\Omega)$  are separable, reflexive,  Banach spaces, (cf. \cite{Rao1}). We recall that  $(\phi_3)$ implies  $(\phi_{3})^{'}$.
\end{rmk}

\nd Using the Poincar\'e inequality \eqref{lambda1} it follows that
\[
\|u\|_\Phi\leq C \|\nabla u\|_\Phi~\mbox{for each}~ u \in W_{0}^{1,\Phi}(\Omega).
\]
holds true for some $C > 0$.
\nd As a consequence,  $\|u\| :=\|\nabla u\|_\Phi$ defines a norm in $W_{0}^{1,\Phi}(\Omega)$, equivalent to $\|.\|_{1,\Phi}$. Let $\Phi_*$ be the inverse of the function
$$
t\in(0,\infty)\mapsto\int_0^t\frac{\Phi^{-1}(s)}{s^{\frac{N+1}{N}}}ds
$$
\nd which extends to ${\mathbb{R}}$ by  $\Phi_*(t)=\Phi_*(-t)$ for  $t\leq 0.$
We say that an N-function $\Psi$ grows  essentially more slowly than $\Phi_*$, we write $\Psi<<\Phi_*$, if
$$
\lim_{t\rightarrow \infty}\frac{\Psi(\lambda t)}{\Phi_*(t)}=0,~~\mbox{for all}~~\lambda >0.
$$

\nd The imbedding below (cf. \cite{A, DT}) will be  used in this paper:
$$
\displaystyle W_{0}^{1,\Phi}(\Omega) \stackrel{\tiny cpt}\hookrightarrow L_\Psi(\Omega)~~\mbox{if}~~\Psi<<\Phi_*,
$$
in particular, as $\Phi<<\Phi_*$ (cf. \cite[Lemma 4.14]{Gz1}),
$$
W_{0}^{1,\Phi}(\Omega) \stackrel{\tiny{cpt}} \hookrightarrow L_\Phi(\Omega).
$$
Furthermore,
$$
W_0^{1,\Phi}(\Omega) \stackrel{\mbox{\tiny cont}}{\hookrightarrow} L_{\Phi_*}(\Omega).
$$

\begin{rmk}\label{rmk-psi}
    The condition $(\psi_1)$ shows that $\Psi<<\Phi_*$, i.e, the function $\Psi$ grows  essentially more slowly than $\Phi_*$.
    In fact, by Proposition  $ \ref{lema_naru}$, stated   below,
    $$\lim_{t\rightarrow\infty}\frac{\Psi(\lambda t)}{\Phi_*(t)}\leq \frac{\lambda^{m_\Psi}}{\Phi_*(1)}\lim_{t\rightarrow\infty}\frac{1}{t^{\ell^*-m_\Psi}}=0,~~\mbox{for all}~~\lambda>0.$$
    In this case $W_{0}^{1,\Phi}(\Omega) \stackrel{cpt}\hookrightarrow L_\Psi(\Omega)$.
\end{rmk}

\nd We refer the reader to   \cite{Fuk_1} for  the  two results below.
\vskip.2cm

\begin{prop}\label{lema_naru}
       Assume that  $\phi$ satisfies  $(\phi_1), (\phi_2),(\phi_3)^{\prime}$.
        Set
         $$
         \zeta_0(t)=\min\{t^\ell,t^m\},~~~ \zeta_1(t)=\max\{t^\ell,t^m\},~~ t\geq 0.
        $$
  \nd Then  $\Phi$ satisfies
       $$
            \zeta_0(t)\Phi(\rho)\leq\Phi(\rho t)\leq \zeta_1(t)\Phi(\rho),~~ \rho, t> 0,
        $$
$$
\zeta_0(\|u\|_{\Phi})\leq\int_\Omega\Phi(u)dx\leq \zeta_1(\|u\|_{\Phi}),~ u\in L_{\Phi}(\Omega).
 $$
    \end{prop}
\begin{prop}\label{lema_naru_*}
    Assume that  $\phi$ satisfies$(\phi_1), (\phi_2),(\phi_3)^{\prime}$.  Set
    $$
    \zeta_2(t)=\min\{t^{\ell^*},t^{m^*}\},~~ \zeta_3(t)=\max\{t^{\ell^*},t^{m^*}\},~~  t\geq 0
    $$
\nd where $1<\ell,m<N$ and $m^* = \frac{mN}{N-m}$, $\ell^* = \frac{\ell N}{N-\ell}$.  Then
        $$
            \ell^*\leq\frac{t^2\Phi'_*(t)}{\Phi_*(t)}\leq m^*,~t>0,
       $$
        $$
            \zeta_2(t)\Phi_*(\rho)\leq\Phi_*(\rho t)\leq \zeta_3(t)\Phi_*(\rho),~~ \rho, t> 0,
       $$
       $$
            \zeta_2(\|u\|_{\Phi_{*}})\leq\int_\Omega\Phi_{*}(u)dx\leq \zeta_3(\|u\|_{\Phi_*}),~ u\in L_{\Phi_*}(\Omega).
        $$
    \end{prop}

\section{ \bf \Large Nehari manifolds  in Orlicz-Sobolev spaces and\\
  Ground State Solutions of problem \eqref{pi} }
\vskip.2cm

\nd  Under the conditions of the present paper   the energy functional J  of \eqref{pi} given by
$$
J(u) = \displaystyle \int_{\Omega} \Phi(|\nabla u|)dx  - \displaystyle \int_{\Omega}F(x, u)dx,~ u \in W^{1, \Phi}_{0}(\Omega),
$$
\nd where
$$
F(x,t) = \int_{0}^{t} f(x,s) ds~\mbox{for}~ s \in \mathbb{R},
$$
\nd  is of class $C^{1}$ and actually
$$
\displaystyle \langle J^{\prime}(u), v \rangle =  \int_\Omega\phi(|\nabla u|)\nabla u\nabla v~dx - \int_\Omega f(x,u) v~dx, ~ u, v\in\w.
$$
\nd Finding weak solutions of  problem \eqref{pi} is equivalent to find critical points of $J$.
\vskip.2cm

\nd The Nehari manifold associated to  $J$ is
given by
$$
\mathcal{N}= \{u \in W^{1,\Phi}_{0}(\Omega) \setminus \{0\}~ |~ \ \langle J^{\prime}(u), u \rangle=0\}.
$$

\nd  The main result of this section is:

\begin{theorem}\label{th1}
Assume  $(\phi_{1})-(\phi_{3})$, $(f_{0}) - (f_{3})$. Then  problem \eqref{pi} admits a nonzero {\rm ground state solution} $u$ in the sense that  $u \in \mathcal{N}$,
$$
  \int_\Omega\phi(|\nabla u|)\nabla u\nabla v~dx=\int_\Omega f(x,u) v~dx, ~v\in\w, ~\mbox{and}
$$
\begin{equation*}
J(u) =   \inf_{w \in \mathcal{N}} J(w) .
\end{equation*}
\end{theorem}

\nd The proof of Theorem \ref{th1}   we will be given  in Section \ref{Pf th1}.
\vskip.3cm

\nd Initially we will establish and prove  a few technical results.
\vskip.2cm

\begin{prop}\label{we}  Assume $(\phi_{1}) -(\phi_{3})$ and $(f_{0})$. Then the functionals
$$
\mbox{\rm{(i)}}~~u \in W^{1, \Phi}_{0}(\Omega) \mapsto \int_\Omega F(x,u)dx, \quad  \mbox{\rm{(ii)}}~~u \in W^{1, \Phi}_{0}(\Omega) \mapsto \int_\Omega f(x,u)udx
$$
\nd are weakly sequentially continuous, w.s.c. for short.
\end{prop}
\nd {\bf Proof of Proposition \ref{we}.}  The  compact embeddings $W^{1, \Phi}_{0}(\Omega) \hookrightarrow L_{\Phi}(\Omega)$ and $W^{1, \Phi}_{0}(\Omega)  \hookrightarrow L_{\Psi}(\Omega)$ will play a crucial role.
Remind  that $\Phi << \Phi_{\star}$ and $\Psi << \Phi_{\star}$.
\vskip.1cm

\nd Let $(u_{n})$ be a sequence in $W^{1, \Phi}_{0}(\Omega)$ such that $u_{n} \rightharpoonup u$ for some $u \in W^{1, \Phi}_{0}(\Omega)$. Hence, up to a subsequence, $u_{n} \rightarrow u$ in $L_{\Psi}(\Omega)$, $u_{n} \rightarrow u$ a.e. in $\Omega$ and $|u_{n}| \leq h$ for some $h \in L_{\Psi}(\Omega)$. Consequently, using $(f_{0})$, we obtain
$$ |f(x,u_{n})u_{n}| \leq C|u_n| + C \Psi(u_{n}) \leq Ch + C \Psi(h) \in L^{1}(\Omega).$$
By the  Dominated Convergence Theorem we have
\begin{equation*}
\lim_{n\rightarrow \infty} \int_{\Omega} f(x, u_{n}) u_{n} dx= \int_{\Omega} f(x,u)u dx.
\end{equation*}
Similarly, one shows  that
\begin{equation*}
\lim_{n\rightarrow \infty} \int_{\Omega} F(x,u_{n})dx = \int_{\Omega} F(x,u)dx.
\end{equation*}
This finishes  the proof.      $\hfill{\rule{2mm}{2mm}}$
\vskip.2cm

\begin{rmk}\label{RMK-1}
\nd  We gather some facts about the behavior of F both near the origin and at infinity. Let $\varepsilon >0$ be a  small number. By $(f_{0})$ and $(f_{2})$  there is a positive constant $C_{\varepsilon}$ such that
$$
|f(x, t)|\leq (\lambda_1-\varepsilon) |t \phi(t)| + C_{\varepsilon}\psi(t), \quad  t \in \Re
$$
and
\begin{eqnarray}\label{aee}
|F(x, t)|\leq (\lambda_1 - \varepsilon) \Phi(t) + C_{\varepsilon} \Psi(t), \quad  t \in \Re.
\end{eqnarray}

\nd In addition, using $ (f_2)$ and $(\phi_3)^{\prime}$ one finds that
\begin{equation}\label{limsup}
\limsup_{t \rightarrow 0} \frac{tf(x, t)}{\Phi(t)} < \frac{\lambda_{1}}{m}.
\end{equation}

\nd By $(f_0)$ and $(\ref{limsup})$ it follows that
\begin{eqnarray}\label{ae-0}
|t f(x, t)|\leq \left(\frac{\lambda_1-\varepsilon}{m}\right) \Phi(t) + C_{\varepsilon}\Psi(t), \quad  t \in \Re.
\end{eqnarray}

\nd As a consequence of  $(\ref{aee})$  and $(\ref{ae-0})$, for each $\varepsilon >0$  there is $C_{\varepsilon}>0$ such that for  each $u \in W^{1, \Phi}_{0}(\Omega)$,
\begin{equation}\label{d2}
\begin{array}{lll}
    \displaystyle\int_\Omega f(x,u)udx & \leq &  \displaystyle\left(\frac{\lambda_1-\varepsilon}{m}\right) \int_{\Omega} \Phi(u)dx +C_\varepsilon \int_{\Omega}  \Psi(u) dx,\\
    \displaystyle\int_\Omega F(x,u)dx & \leq & \displaystyle(\lambda_1 - \varepsilon) \int_{\Omega} \Phi(u)dx +C_\varepsilon \int_{\Omega} \Psi(u) dx.
\end{array}
\end{equation}

\vskip.1cm

\nd Using the embedding $W^{1,\Phi}_{0}(\Omega)  \hookrightarrow  L_{\Psi}(\Omega) $ and Proposition $\ref{lema_naru}$ it follows that
\begin{equation*} \label{d22}
\int_\Omega f(x,u)u dx\leq \displaystyle\left(\frac{\lambda_1-\varepsilon}{m}\right) \int_{\Omega} \Phi(u)dx + C_\varepsilon \max ( \|u\|^{\ell_{\Psi}}, \|u\|^{m_{\Psi}} )
\end{equation*}
and
\begin{equation}\label{d222}
 \int_\Omega F(x,u) dx \leq \displaystyle(\lambda_1- \varepsilon) \int_{\Omega} \Phi(u)dx + C_{\varepsilon} \max ( \|u\|^{\ell_{\Psi}}, \|u\|^{m_{\Psi}} )
\end{equation}
\nd for  $u \in W^{1, \Phi}_{0}(\Omega)$.
\end{rmk}
\vskip.3cm

\nd At this point, aiming to determine   the behavior of $J$ on  $\mathcal{N}$ we introduce the fibering maps $\gamma_{u} : (0, \infty) \rightarrow \mathbb{R}$ for $u \in \w\setminus \{ 0\}$, defined by
\begin{equation*}
\gamma_{u}(t) = J (t u),~~ t \in (0,\infty),
\end{equation*}
(see \cite{Brown2, Brown1}).  In this regard we shall study the behavior of $\gamma_{u}(t)$  for  both  $t$ near infinity and $t$ near the origin.
\vskip.2cm

\nd Our next result is:

\begin{prop}\label{be}
Assume $(\phi_{1}) -(\phi_{3})$,~~ $(f_{0}) -(f_{3})$ and let $u \in \w\setminus \{ 0\}$. Then
\begin{equation*}
\mbox{\rm{(i)}}~~ \lim_{t \rightarrow 0} \dfrac{\gamma_{u}(t)}{t^{m}} > 0,~~~~ \lim_{t \rightarrow \infty} \dfrac{\gamma_{u}(t)}{t^m} = - \infty,
\end{equation*}
and in addition,
\begin{equation*}
\mbox{\rm{(ii)}}~~ \lim_{t \rightarrow 0} \dfrac{\gamma^{\prime}_{u}(t)}{t^{m-1}} >  0,~~~~ \lim_{t \rightarrow \infty} \dfrac{\gamma^{\prime}_{u}(t)}{t^{m-1}} = - \infty.
\end{equation*}
\end{prop}

\nd {\bf Proof of Proposition \ref{be}.}  By \eqref{d222} one infers  that
\begin{equation*}
\gamma_{u}(t) \geq \int_{\Omega} \Phi(|\nabla (t u)|) dx - (\lambda_1 - \varepsilon) \int_{\Omega} \Phi(tu) dx - C_{\varepsilon} \max \left(\|tu\|^{\ell_{\Psi}} ,\|(t u)\|^{m_{\Psi}}\right).
\end{equation*}
Applying the Poincar\'e inequality we have
\begin{equation*}
\gamma_{u}(t) \geq \left( 1 - \dfrac{\lambda_1 - \varepsilon}{\lambda_{1}} \right) \int_{\Omega} \Phi(|\nabla t u|) dx  - C_{\varepsilon} \max \left(\|tu\|^{\ell_{\Psi}} ,\|tu\|^{m_{\Psi}}\right).
\end{equation*}

\nd Appying Proposition \ref{lema_naru} it follows that
\begin{equation*}
\gamma_{u}(t) \geq \left( 1 - \dfrac{\lambda_1 - \varepsilon}{\lambda_{1}} \right) t^{m} \int_{\Omega} \Phi(|\nabla u|) dx  - C_{\varepsilon} \max \left(\|tu\|^{\ell_{\Psi}}, \|tu\|^{m_{\Psi}}\right).
\end{equation*}

\nd By the arguments above we have: for each  $u \in \w\backslash \{ 0\}$ and   $0 < t < 1$,
\begin{equation}\label{limm}
\dfrac{\gamma_{u}(t)}{t^{m}} \geq \left( 1 - \dfrac{\lambda_1 - \varepsilon}{\lambda_{1}} \right) \int_{\Omega} \Phi(|\nabla u|) dx  - C_{\varepsilon} \dfrac{\max \left(\|tu\|^{\ell^{\Psi}}, \|tu\|^{m_{\Psi}}\right)}{t^{m}}.
\end{equation}
Using the fact that $m < \ell_{\Psi}$,    \eqref{limm}  rewrites as
\begin{equation*}
\dfrac{\gamma_{u}(t)}{t^{m}} \geq \left( 1 - \dfrac{\lambda_1 - \varepsilon}{\lambda_{1}} \right) \int_{\Omega} \Phi(|\nabla u|) dx + o(1)
\end{equation*}
where $o(1)$ denotes a quantity that goes to zero as $t \rightarrow 0$.

\nd Hence
\begin{equation*}
\lim_{t \rightarrow 0}\dfrac{\gamma_{u}(t)}{t^{m}} \geq \left( 1 - \dfrac{\lambda_1 - \varepsilon}{\lambda_{1}} \right) \int_{\Omega} \Phi(|\nabla u|) dx > 0.
\end{equation*}

\nd Next  we shall compute the limit of $\gamma_{u}(t)/t^{m}$ at infinity. Using Proposition \ref{lema_naru} we get
\begin{equation*}
\frac{\gamma_{u}(t)}{t^{m}} \leq \int_{\Omega} \Phi(|\nabla  u|) dx -  \dfrac{1}{t^{m}} \int_{\Omega} F(x, t u)dx.
\end{equation*}
\nd Applying Fatou's Lemma and  $(f_{3})$ we  get
\begin{equation*}
\lim_{t \rightarrow \infty} \dfrac{\gamma_{u}(t)}{t^{m}} \leq  \int_{\Omega} \Phi(|\nabla u|) dx - \liminf_{t \rightarrow \infty} \int_{\Omega} \dfrac{F(x, tu)}{t^{m}} dx = - \infty.
\end{equation*}
\nd We emphasize that $(f_{0})$ and $(f_{3})$ ensure that
\begin{equation*}
\liminf_{t \rightarrow \infty} \int_{\Omega} \dfrac{F(x, tu)}{t^{m}} dx = \liminf_{t \rightarrow \infty} \int_{\Omega} \dfrac{F(x, tu)}{|t u|^{m}}|u|^{m} dx  = + \infty.
\end{equation*}

\nd The limits involving $\gamma_{u}^{\prime}(t)/t^{m -1}$ are computed by arguments similar to the ones above. This ends the proof.  $\hfill{\rule{2mm}{2mm}}$

\begin{prop}\label{lll}
Suppose $(\phi_{1}), (\phi_{2})$, $(\phi_{3})$. Then the functionals
$$
\mbox{\rm (i)}~~~u  \in \w \mapsto \int_{\Omega} \phi(|\nabla u|)|\nabla u|^{2} dx
$$
\nd and
$$
 \mbox{\rm (ii)}~~ u \in \w \mapsto  \int_{\Omega} (m \Phi(|\nabla u|) - \phi(|\nabla u|)|\nabla u|^{2} ) dx
$$
are weakly sequentially lower semicontinuous,  w.s.l.s.c. for short.
\end{prop}

\nd {\bf Proof of Proposition \ref{lll}.}  It is an easy matter to show that both functionals are  lower semicontinuous. So it is enough to show that they are also convex. To this end consider the functions
$$
L_{1}(t) = \phi(t)t^{2},~~~~~~   L_{2}(t) = m \Phi(t)-\phi(t)t^{2},~~~ t \geq 0.
$$
It follows by an easy computation, using  condition $(\phi_{3})$ that
$$
L^{\prime \prime}_{1}(t) \geq 2 (t\phi(t))^{\prime} + (\ell -2) (t\phi(t))^{\prime} = l(t\phi(t))^{\prime} \geq 0 .
$$
On the other hand,  this time using condition  $(\phi_{3})^{\prime}$ we infer that
$$
L^{\prime \prime}_{2}(t) \geq (m - \ell)(t \phi(t))^{\prime} \geq 0.
$$
Thus  $L_{1}, L_{2}$ are convex functions. As a consequence the functionals in (i)-(ii)  above are  convex. This finishes the proof.    $\hfill{\rule{2mm}{2mm}}$
\vskip.2cm

\nd Next we shall discuss the relation between the fibering map $\gamma_{u}$ and  $\mathcal{N}$. Roughly speaking it will be shown that  the map  $\gamma_{u}$ crosses  $\mathcal{N}$  once in a suitable sense.

\begin{prop}\label{l1}
Assume $(\phi_{1}) -(\phi_{3})$,  $(f_{0})- (f_{3})$. Then for each $u \in W^{1,\Phi}_{0}(\Omega) \setminus \{0\}$ there is  only a  $t = t(u) >0$ such that $tu \in \mathcal{N}$. Moreover, $J(u)> 0$ for each $u \in \mathcal{N}$.
\end{prop}

\nd {\bf Proof of Proposition \ref{l1}.}   Let $u \in W^{1,\Phi}_{0}(\Omega) \setminus \{0\}$. At first, note that  by the very  definition of  $\gamma_{u}(t)$, $tu \in \mathcal{N}$ if and only if $\gamma_{u}'(t)=0$.
\vskip.1cm

\nd On the other hand,  by Proposition \ref{be} we have
$$
\gamma^{\prime}_{u}(t) > 0~~~\mbox{for}~~~  t~~  \mbox{small enough}
$$
\nd and
$$
\gamma^{\prime}_{u}(t) < 0~~~ \mbox{for}~~~  t~~~ \mbox{big enough}.
$$

\nd Since the map  $t \mapsto \gamma^{\prime}_{u}(t)$ is  continuous  there is at least one  number $t \in (0,\infty)$  such  that $\gamma^{\prime}_{u}(t) = 0$. This means that $t u \in \mathcal{N}$.
\vskip.1cm

\nd We claim that  there is only one  $t = t(u)$ such that $\gamma_{u}^{\prime}(t) = 0$. Indeed,  recall that $\gamma_{u}^{\prime}(t) = \langle J^{\prime}(tu), u \rangle$. So
\begin{equation}\label{e1}
\dfrac{d}{d t} \left[\dfrac{\gamma^{\prime}_u(t)} {t^{m -1}} \right]= \int_{\Omega} \dfrac{d}{d t} \left[\dfrac{\phi(|\nabla t u|)\nabla t u \nabla u}{t^{m-1}}\right]  dx
- \int_{\Omega}\dfrac{d}{dt} \left[\dfrac{f(x, t u) u}{t^{m-1}}\right] dx.
\end{equation}

\nd At this point  we remark  that $(\phi_{3})$ implies
\begin{equation*}
\ell -2 \leq \inf_{t > 0} \dfrac{t \phi^{\prime}(t)}{\phi(t)} \leq \sup_{t > 0} \dfrac{t \phi^{\prime}(t)}{\phi(t)} \leq m -2.
\end{equation*}

\nd Using the inequalities just above we infer that  for each $x \in \Omega$ and  $t > 0$,
\begin{equation}\label{e2}
\dfrac{d}{d t} \left[\dfrac{\phi(|\nabla t u|)\nabla t u \nabla u}{t^{m-1}}\right] = \dfrac{|\nabla u|^{2}\left[ \phi^{\prime}(|\nabla t u|)|\nabla t u|
-(m-2) \phi(|\nabla t u|) \right]}{t^{m -1}} \leq 0.
\end{equation}

\nd It follows by using \eqref{e1} and \eqref{e2} that
\begin{equation}\label{e3}
\dfrac{d}{d t} \left[\dfrac{\gamma^{\prime}_u(t)} {t^{m -1}}\right] \leq - \int_{\Omega} \dfrac{d}{dt} \left[\dfrac{f(x, t u) u}{t^{m-1}}\right]  dx
= - \int_{\Omega} \dfrac{d}{dt} \left[\dfrac{f(x, t u)}{|tu|^{m-2} tu}  \right] |u|^{m}  dx .
\end{equation}

\nd Now, using $(f_{1})$ and \eqref{e3} we get
\begin{equation*}
\dfrac{d}{d t} \left[\dfrac{\gamma^{\prime}_u(t)} {t^{m -1}}\right] < 0
\end{equation*}
\nd for each  $t > 0$ and $u \in \w\backslash \{0\}$.
\vskip.1cm

\nd Therefore  $t \mapsto \dfrac{\gamma^{\prime}_u(t)} {t^{m-1}}$ is
a decreasing  function that vanishes  once in $(0,\infty)$  so that  there is an only  $t = t(u) > 0$  such  that
${\gamma^{\prime}_u(t)}/ {t^{m-1}}= 0$.
\vskip.1cm

\nd  Thus the function $\gamma_{u}$ admits a unique critical point namely $t = t(u) >  0$ and actually,  $t u \in \mathcal{N}$. Moreover it follows by Proposition \ref{be}  that $t(u)$ is a  maximum point of $\gamma_u$ on $(0,\infty)$ and, in fact $\gamma_u(t(u))>0$, which implies that  $J(t(u)u)>0$.  The arguments above also show that  $\gamma_{u}^{\prime \prime}(t) < 0$ for each  $u \in \w \setminus \{0\}$.
\vskip.1cm

\nd   Finally,  since $u \in \mathcal{N}$ if only if   $t(u)=1$, we deduce that $J(u)>0$ for each $u \in \mathcal{N}$. This completes the proof.       $\hfill{\rule{2mm}{2mm}}$
\vskip.2cm

\begin{prop} \label{deri}
\nd  Assume $(\phi_{1})-(\phi_{3})$, $(f_{0}) - (f_{3})$. Then  $I:\w\to \mathbf{R}$,
$$I(u)=\int_\Omega \phi(|\nabla  u|)|\nabla u|^2dx,~~u\in\w ,$$
is  ${C}^1$ and
$$\langle I'(u),v\rangle=\int_\Omega [2\phi(|\nabla u|)+\phi'(|\nabla u|)|\nabla u|]\nabla u.\nabla v dx,~u,v\in\w .
$$
\end{prop}

\nd {\bf Proof of Proposition  \ref{deri}.}  Set
$$
g(t)=\phi(|\nabla u+t\nabla v|)|\nabla u+t\nabla v|^2,~~0 \leq t \leq 1.
$$
\nd It follows that  $g\in {C}^1$ and, actually
 $$
g'(t)=[2\phi(|\nabla u+t\nabla v|)+\phi'(|\nabla u+t\nabla v|)|\nabla u+t\nabla v|](\nabla u+t\nabla v).\nabla v.
$$

\nd  In addition,  there is  $\theta \in \mathbf{R}$ with $0< \theta<t \leq 1$  such that
    $$
\frac{g(t)-g(0)}{t}=g'(\theta).
$$
\nd Thus
  \begin{eqnarray*}
        \langle I'(u),v\rangle&=&\lim_{t\rightarrow0}\int_\Omega \frac{\phi(|\nabla u+t\nabla v|)|\nabla u+t\nabla v|^2-\phi(|\nabla u|)|\nabla u|^2}{t}dx\\
      & = & \lim_{\theta \rightarrow0}\int_{\Omega} g'(\theta) dx.
       \end{eqnarray*}

\nd {\bf Claim.}   There is  $h \in L^{1}(\Omega)$ such  that $|g'(\theta)| \leq h$.
\vskip.1cm

\nd  At first we recall that by  $(\phi_{3})^{\prime}$,
     \begin{equation}\label{limitephilinha-0}
        |\phi'(t)t|\leq \max\{|\ell-2|,|m-2|\}\phi(t),~~  0 \leq t < \infty
    \end{equation}
 \nd and  as a consequence,
    \begin{equation}\label{limitephilinha}
        \lim_{t\rightarrow0}|\phi'(t)t^2|=\lim_{t\rightarrow 0} [\max\{|\ell-2|,|m-2|\}\phi(t) t] = 0.
    \end{equation}

\nd Using  $(\phi_{1}),(\phi_{2}),(\phi_{3})$ and $0 \leq \theta \leq 1$ we have,
\begin{eqnarray*}
       \left|g'(\theta)\right| & \leq & [2\phi(|\nabla u+\theta \nabla v|)+|\phi'(|\nabla u+\theta \nabla v|)||\nabla u+\theta \nabla v|]|\nabla u+\theta \nabla v||\nabla v|\\
       & \leq & [2+\max(|\ell -2|,|m-2|)]\phi(|\nabla u|+|\nabla v|)(|\nabla u|+|\nabla v|)|\nabla v| .
    \end{eqnarray*}
    \vskip.2cm

\nd  Next we show that  $[2+\max(|\ell -2|,|m-2|)]\phi(|\nabla u|+|\nabla v|)(|\nabla u|+|\nabla v|)|\nabla v| \in  L^1(\Omega$).
\vskip.3cm

\nd Indeed,  using Young's inequality, the inequality $\widetilde \Phi(t\phi(t))\leq \Phi(2t)$ and the fact that  $\Phi\in \Delta_2$ we have
    \begin{eqnarray*}
        \phi(|\nabla u|+|\nabla v|)(|\nabla u|+|\nabla v|)|\nabla v| & \leq & \Phi(|\nabla v|)+\widetilde \Phi(\phi(|\nabla u|+|\nabla v|)(|\nabla u|+|\nabla v|))\\
        & \leq & \Phi(|\nabla v|)+ \Phi(2(|\nabla u|+|\nabla v|))\\
        & \leq & \Phi(|\nabla v|)+ 2^m\Phi(|\nabla u|+|\nabla v|)\\
        & \leq & (1+2^m)\Phi(|\nabla u|+|\nabla v|) \in  L^1(\Omega).
    \end{eqnarray*}

\nd So,
$$
|g'(\theta)|  \leq  | [2+\max(|\ell -2|,|m-2|)](1+2^m)\Phi(|\nabla u|+|\nabla v|)|] := h \in L^1(\Omega).
$$

\nd This ends the proof of the claim.

\nd    It remains to show that $I'$ is continuous. Indeed, let $(u_n)\subseteq \w$ such that $u_n\rightarrow u$ em $\w$. Then
   \begin{eqnarray*}
     \int_\Omega \Phi(|\nabla u_n- \nabla u|)dx  \stackrel{n\rightarrow\infty} \longrightarrow  0,\\
\\
\nabla u_n\rightarrow \nabla u~~\mbox{ a.e in }~~ \Omega,\\
 \\
 |\nabla u_n|\leq h_1~ \mbox{a.e. in}~ \Omega~ \mbox{for some}~  h_1\in L^1(\Omega).
\end{eqnarray*}

\nd By  arguments  as above,
    \begin{eqnarray*}
         |[2\phi(|\nabla u_n|)+\phi'(|\nabla u_n|)|\nabla u_n|]\nabla u_n.\nabla v|
& \leq & [2\phi(|\nabla u_n|)+|\phi'(|\nabla u_n|)||\nabla u_n|]|\nabla u_n||\nabla v|\\
         & \leq & [2\phi(|\nabla u_n|)|\nabla u_n|+|\phi'(|\nabla u_n|)||\nabla u_n|^2]|\nabla v|\\
         & \leq &C \phi(|\nabla u_n|)|\nabla u_n||\nabla v|\\
         & \leq & C]\phi(h_1)h_1|\nabla v|\in L^1(\Omega),
    \end{eqnarray*}
\nd    where $C:= 2+\max(|\ell-2|,|m-2|)$.
\vskip.2cm

\nd  Since $\nabla u_n\rightarrow \nabla u$ a.e. in $\Omega$ we get by (\ref{limitephilinha-0})-(\ref{limitephilinha}) that,
 $$
[2\phi(|\nabla u_n|)+\phi'(|\nabla u_n|)|\nabla u_n|]\nabla u_n.\nabla v\rightarrow[2\phi(|\nabla u|)+\phi'(|\nabla u|)|\nabla u|]\nabla u.\nabla v~\mbox{a.e. in}~ \Omega
$$

\nd Applying the Lebesgue dominated convergence theorem,
    $$\lim_{n\rightarrow +\infty}\langle I'(u_n)-I'(u),v\rangle=0.$$
    Hence we have $I \in {C}^1(\w;\mathbb{R})$ which completes the proof of Proposition \ref{deri}. $\hfill{\rule{2mm}{2mm}}$

\begin{prop}\label{c11}
Assume  $(\phi_{1})- (\phi_{3})$, $(f_{0})- (f_{3})$. Then  $\mathcal{N}$ is a ${C}^1$-submanifold of $W^{1,\Phi}_{0}(\Omega)$. In addition,  any critical point of $J_{\mid \mathcal{N}}$ is a critical point of $J$.
\end{prop}

\nd {\bf Proof of Proposition  \ref{c11}.} Note that by the very definition of $\gamma_u$,
$$
\gamma_u'(t)= \langle J'(tu), u \rangle,~ u \in  W^{1,\Phi}_{0}(\Omega) \backslash \{ 0\}.
$$
\nd Consider  the functional  $J_{t}: \w(\Omega) \rightarrow \mathbb{R}$ defined  by
\begin{equation*}
J_{t}(u) = I_{t}(u) - \int_{\Omega} f(x, t u) u dx,
\end{equation*}
\nd where
\begin{equation*}
I_{t}(u) = \int_{\Omega} \phi(|\nabla (t u)|) \nabla (t u) \nabla u dx.
\end{equation*}
\nd Using arguments as in the proof of  Proposition \ref{deri} one shows that  $I_{t} \in C^{1}$ and
\begin{equation*}
\langle I^{\prime}_{t}(u), v \rangle = \int_{\Omega} \left[ 2 \phi(|\nabla (t u)|) + \phi^{\prime}(|\nabla (t u)|) |\nabla (t u)| \right] \nabla u \nabla v dx,~~ u, v \in \w(\Omega),~t \in \mathbb{R}.
\end{equation*}
\nd One also shows  that
\begin{equation*}
\gamma_u''(t)= \int_{\Omega} \left[ \phi(|\nabla (t u)|) + \phi^{\prime}(|\nabla (t u)|) |\nabla (t u)| \right] |\nabla u|^{2} - \int_{\Omega} f^{\prime}(x, t u) u^{2} dx.
\end{equation*}
\nd Set
$$
R(u) = \langle J^{\prime}(u), u \rangle,~u \in \w(\Omega).
$$
\nd It follows that  $R \in C^{1}$, (see Proposition \ref{deri}). Actually  since  $t = 1$ is the global maximum of $\gamma_{u}$, see Proposition \ref{be} $(i)$ and
Proposition \ref{l1}, we observe that
$$
\langle R'(u),u\rangle = \gamma_u''(1) < 0,~  u \in \mathcal{N}.
$$
\nd  Using the fact that $\mathcal{N} = R^{-1}(0)$ and $0$ is a regular value for $R$, the set $\mathcal{N}$ is a $C^1$-submanifold of $W^{1,\Phi}_{0}(\Omega)$.
\vskip.1cm

\nd To finish the proof, we assume that $u \in \mathcal{N}$ is a critical point of $J_{\mid\mathcal{N}}$.  Applying the Lagrange Multiplier Theorem, we have
 $$
J'(u)= \mu R'(u)~\mbox{ for some}~ \mu \in \Re.
$$
\nd  Taking $u$ as a test function it follows that
$$
\mu \langle R'(u),u\rangle =\langle J'(u),u\rangle=0.
$$
\nd Reminding that    $\langle R'(u),u\rangle = \gamma_u''(1) < 0~  \mbox{for}~ u \in \mathcal{N}$ we infer that  $\mu=0$. Therefore $J'(u) \equiv 0$, so that  $u$ is a free critical point of $J$. This completes the proof.   $\hfill{\rule{2mm}{2mm}}$

\begin{prop}\label{l22}
Assume  $(\phi_{1})-(\phi_{3})$, $(f_{0})- (f_{3})$. Then there is a constant $C>0$ such that $\|u\| \geq C$ for each $u \in \mathcal{N}$.
\end{prop}

\nd {\bf Proof of Proposition  \ref{l22}.} Assuming the contrary, there is  $(u_{n}) \subset  \mathcal{N}$  such  that
$\|u_{n}\| \leq \dfrac{1}{n}$ for each integer $n \geq 1$.
Let $\epsilon > 0$. Using \eqref{d2} and the Poincar\'e inequality we find  some $C_{\epsilon} > 0$ such that
\begin{eqnarray*}
\int_{\Omega} \Phi(|\nabla u_{n}|) dx &\leq& \dfrac{1}{\ell} \int_{\Omega} \phi(|\nabla u_{n}|) |\nabla u_{n}|^{2}dx  \nonumber \\
&=&  \dfrac{1}{\ell} \displaystyle\int_{\Omega}f(x,u_{n})u_{n}dx \leq \dfrac{\lambda_1 - \varepsilon}{\lambda_{1}} \int_{\Omega} \Phi(|\nabla u_{n}|) dx + C_{\varepsilon} \int_{\Omega} \Psi(u)dx \nonumber \\
\end{eqnarray*}
\nd Hence
\begin{equation*}
\left(1 - \dfrac{\lambda_1 - \varepsilon}{\lambda_{1}} \right) \int_{\Omega} \Phi(|\nabla u_{n}|) dx \leq C_{\varepsilon} \int_{\Omega} \Psi(u)dx
\end{equation*}
\nd so that
\begin{equation*}
\int_{\Omega} \Phi(|\nabla u_{n}|) dx \leq \left(1 - \dfrac{\lambda_1 - \varepsilon}{\lambda_{1}} \right)^{-1} C_{\varepsilon} \int_{\Omega} \Psi(u)dx.
\end{equation*}
\nd Applying  Proposition \ref{lema_naru} we find
\begin{equation*}
 \| u_{n}\|^{m} \leq \int_{\Omega} \Phi(|\nabla u_{n}|) dx \leq C_{\varepsilon} \max \left(\|u_{n}\|^{m_{\Psi}}, \|u_{n}\|^{\ell_{\Psi}} \right) = C_{\varepsilon}  \|u_{n}\|^{\ell_{\Psi}}.
\end{equation*}
\nd  Dividing the last expression by $\| u_{n}\|^{m}$ we get to
\begin{equation*}
1 \leq C_{\varepsilon} \|u_{n}\|^{\ell_{\Psi} - m}.
\end{equation*}
\nd Passing to the limit  as $n \rightarrow \infty$ and using the fact that  $\ell_{\Psi} > m$  we get to a contradiction. So the Nehari manifold $\mathcal{N}$ is bounded away from zero by some positive constant $C$. This ends the proof.  $\hfill{\rule{2mm}{2mm}}$

\section  { \bf \Large Proof of Theorem \ref{th1}}\label{Pf th1}

\nd At first we will establish a few Lemmas. Set
$$
c_{\mathcal{N}}:=\inf_{\mathcal{N}} J.
$$
\nd The first lemma  establishes that any minimizing sequence  is bounded in $\w$.

\begin{lem}\label{l33} Assume $(\phi_{1})-(\phi_{3})$, $(f_{0})- (f_{3})$.
Let $(u_{n}) \subset  \mathcal{N}$ be a minimizing sequence of $J$ over the Nehari manifold ${\mathcal{N}}$, that is,    $(u_{n}) \subset {\mathcal{N}}$ satisfies $J(u_n) \rightarrow c_{\mathcal{N}}$. Then $(u_n)$ is  bounded  in $W^{1,\Phi}_{0}(\Omega)$.
\end{lem}

\nd {\bf Proof of Lemma \ref{l33}.} Let $(u_n) \subset \mathcal{N}$ be a minimizing  sequence, that is,  $(u_{n}) \subset {\mathcal{N}}$ and  $J(u_n) \rightarrow c_{\mathcal{N}}$. Assume on the contrary that   $\|u_n\|\rightarrow \infty$.

\nd Set  $v_{n}=\frac{u_{n}}{\|u_{n}\|}$. Then  $\|v_{n}\| = 1$  and there is $v \in W^{1,\Phi}_{0}(\Omega) $ such that $v_n \rightharpoonup v$ in $W^{1,\Phi}_{0}(\Omega)$.

\nd We claim that $v \neq 0$. Assume, by the way of contradiction,  that  $v \equiv 0$. Since  $u_{n} \in {\mathcal{N}}$ it follows that
$$
J(u_n)= \max_{t>0} J(tu_n)~ \mbox{ for each}~ n.
$$
\nd  Let $M > 0$ be a constant.  Now we observe that
$$
c_{\mathcal{N}}+o_{n}(1)= J(u_{n})\geq J(M v_{n})
=  \int_{\Omega} \Phi(|\nabla (M v_{n})|) dx - \displaystyle \int_{\Omega} F(x, Mv_{n})dx.
$$
\nd Since $v_n \rightharpoonup  0$ in $W^{1,\Phi}_{0}(\Omega)$ it follows by Proposition \ref{we} (i) that
$$
\displaystyle\int_{\Omega}F(x,M v_{n}) dx \rightarrow 0.
$$
\nd Employing  Proposition \ref{lema_naru} we have
$$
c_{\mathcal{N}}+o_{n}(1)\geq  J(u_{n}) = \int_{\Omega} \Phi(|\nabla (M v_{n})|) dx + o_{n}(1) \geq  \min(M^{\ell}, M^{m}) + o_{n}(1),~  M>0.
$$
\nd Passing to the limit in the inequalitities just above we get
$$
c_{\mathcal{N}}  \geq  \min(M^{\ell}, M^{m}),~ M>0,
$$
\nd which is impossible.  Therefore  $v \neq 0$.

\nd Remember we are assuming that  $\|u_n\|\rightarrow \infty$ and $J(u_n) \rightarrow c_{\mathcal{N}}$. Hence
$$
\frac{J(u_{n})}{\|u_{n}\|^{m}} =
o_{n}(1).
$$
\nd Applying  Proposition \ref{lema_naru} we have
\begin{eqnarray*}
\displaystyle \int_{\Omega}\frac{F(x,u_{n})}{\|u_{n}\|^{m}} dx
&=&\dfrac{1}{\|u_{n}\|^{m}}\int_{\Omega} \Phi(|\nabla u_{n}|) dx  +o_{n}(1)\\
&\leq&  \int_{\Omega} \Phi(|\nabla v_{n}|) dx +o_{n}(1) =1 + o_{n}(1)
\end{eqnarray*}
\nd Passing to the limit above we have
\begin{equation*}\label{ee1}
\limsup_{n \rightarrow \infty} \displaystyle \int_{\Omega}\frac{F(x,u_{n})}{\|u_{n}\|^{m}}dx \leq 1
\end{equation*}
\nd On the other hand, it follows by  $(f_{3})$  and  L'Hospital rule that
$$
 \displaystyle \lim_{t \to \infty}\frac{F(x, t)}{t^{m}} = + \infty.
$$
\nd Aplying  Fatou's Lemma and using the fact that  $v \not \equiv 0$, we have
\begin{eqnarray*}
\liminf_{n \rightarrow \infty} \int_{\Omega} \dfrac{F(x,u_{n})}{\|u_{n}\|^{m}} dx &\geq& \int_{\Omega} \liminf_{n \rightarrow \infty} \left\{\dfrac{F(x,u_{n})}{\|u_{n}\|^{m}}\right\} dx \nonumber \\
&=&  \int_{\Omega} \liminf_{n \rightarrow \infty} \left\{\dfrac{F(x, u_{n})}{|u_{n}|^{m}} |v_{n}|^{m} \right\}dx = + \infty,  \nonumber \\
\end{eqnarray*}
\nd which is impossible. Thus  $(u_{n})$ is  bounded in $W^{1,\Phi}_{0}(\Omega)$. The proof is  complete. $\hfill{\rule{2mm}{2mm}}$

\begin{lem}\label{pr1}
\label{l4}
Assume $(\phi_{1})-(\phi_{3})$, $(f_{0})- (f_{3})$. Then there exists $u \in \mathcal{N}$ such that $$
c_{\mathcal{N}} = J(u) > 0.
$$
\end{lem}
\nd {\bf Proof of Lemma \ref{pr1}.}~  Let $(u_{n}) \subset \mathcal{N}$ be a minimizing sequence for $J$ over $\mathcal{N}$.
By  Lemma \ref{l33}, there is $u \in \w$ such that
$$
u_{n}\rightharpoonup u \ \ \mbox{in} \ \ W^{1,\Phi}_{0}(\Omega) .
$$
\nd {\bf Claim.}      $u \not \equiv 0$.
\vskip.1cm

\nd Indeed, assume on the contrary that  $u \equiv 0$.  Using condition   $(\phi_3)$  we obtain
$$
(\ell - 2)  (s \phi(s))^{\prime} \leq (s \phi(s))^{\prime \prime} s \leq (m-2) (s \phi(s))^{\prime}.
$$
\nd Integrating  from $0$ to $t$ in the inequalities above, term by term, we get to
$$
(\ell - 2) t \phi(t)   \leq (t \phi(t))^{\prime} t - t \phi(t) \leq (m-2) t \phi(t).
$$
\nd Now we get
$$
(\ell - 1) t \phi(t) \leq (t \phi(t))^{\prime} t \leq (m-1) t \phi(t),
$$
\nd which gives
$$
(\ell - 1) t \phi(t) \leq {(t \phi(t))^{\prime} t}\leq (m-1) t \phi(t).
$$
\nd Applying arguments like in \cite{kaye} we get to
$$
\ell \leq \dfrac{t^2 \phi(t)}{\Phi(t)} \leq m,
$$
\nd which gives
$$
\Phi(t) \leq \dfrac{1}{\ell}{t^2 \phi(t)}.
$$
\nd Using the inequality just above, the fact that $u_n \in \mathcal{N}$  we have
\begin{equation}\label{P3.2}
0 \leq \int_{\Omega} \Phi(|\nabla u_{n}|) dx \leq \dfrac{1}{\ell} \int_{\Omega} \phi(|\nabla u_{n}|) |\nabla u_{n}|^{2}dx = \dfrac{1}{\ell}\int_{\Omega} f(x,u_{n}) u_{n}.
\end{equation}
\nd Applying  Proposition \ref{we} (ii),  we get
$$
\int_{\Omega} f(x,u_{n}) u_{n} dx = o_{n}(1).
$$

\nd As a consequence of (\ref{P3.2}),  $\|u_{n}\|\rightarrow 0$,  contradicting  Proposition \ref{l22}. Therefore  $u \not \equiv 0$, proving the {\bf Claim}.
\vskip.1cm

\nd  As a consequence of Propositions  \ref{lll}, \ref{we} (ii) and (\ref{lll}) (i), we have
$$
u \in W^{1,\Phi}_{0}(\Omega) \mapsto \langle J'(u), u \rangle~~\mbox{is w.s.l.s.c.}.
$$
\nd Hence
$$
\langle J'(u),u \rangle  \leq \liminf_{n\rightarrow \infty} \langle J'(u_{n}),u_{n} \rangle =0.
$$
\nd Recall  that $\gamma_u^{\prime}(1) = \langle J'(u),u \rangle \leq 0$. By Proposition \ref{l1}  and its proof there is $t \in(0,1]$ such  that $\gamma_u^{\prime}(tu)  = 0$. Hence    $tu  \in \mathcal{N}$.
\vskip.1cm

\nd We claim that  $t = 1$ so that  $u$ is in $\mathcal{N}$.
\vskip.1cm

\nd Indeed,  assume on the contrary, that $t \in(0,1)$. In this case we get
\begin{eqnarray}\label{es}
c_{\mathcal{N}}&\leq& J(tu)= J(tu)- \frac{1}{m} \langle J'(tu ),tu\rangle \nonumber \\
&=& \displaystyle \int_{\Omega} \Phi(|\nabla (t u)|)- \frac{1}{m} \phi(|\nabla (t u)|) |\nabla (t u)|^{2}   dx  + \displaystyle \int_{\Omega}\left\{\frac{1}{m} f(x,tu)tu -
F(x, tu)\right\}. \nonumber \\
\end{eqnarray}

\nd Using $(f_{1})$, we get
$$
f'(x,t)t - \frac{1}{m} f(x,t) > 0,~  t>0.
$$
\nd  But the inequality above  implies that
\begin{eqnarray}\label{esti}
t \mapsto \frac{1}{m} f(x,t)t- F(x, t) \ \ \mbox{is increasing in }  (0,\infty) \, \, \mbox{for each} \, \, x \in \Omega.
\end{eqnarray}
Indeed,  using $(f_1)$  we have
\begin{equation*}
\dfrac{d}{d t} \left\{ \frac{1}{m} f(x,t)t- F(x ,t) \right\} =  t^{m} \dfrac{d}{d t} \left\{ \dfrac{f(x, t)}{t^{m -1}} \right\}, t > 0, x \in \Omega.
\end{equation*}

\nd We also have that
$$
t \mapsto \Phi(|\nabla (t u)|)- \frac{1}{m} \phi(|\nabla (t u)|) |\nabla (t u)|^{2}~~\mbox{ is increasing on}~ (0, \infty).
$$
\nd Indeed, setting
 $$
L_{1}(t) = m \Phi(t) - t^{2} \phi(t),~ t > 0
$$
\nd  we find that
$$
L_{1}^{\prime}(t) = (m -1) t \phi (t) - t (t \phi(t))^{\prime}.
$$
Now we observe that by $(\phi_{3})$,
 $$
(\ell - 1) \phi(t) \leq  (t \phi(t))^{\prime} \leq (m -1) \phi(t),~t > 0.
$$
\nd  This shows that $L_{1}$ is increasing.

\nd At this point,  using $\eqref{es}$ and $\eqref{esti}$ we conclude that
$$
c_{\mathcal{N}} < \displaystyle \int_{\Omega} \left\{ \Phi(|\nabla u|)- \frac{1}{m} \phi(|\nabla u|) |\nabla u|^{2} \right\}dx + \displaystyle \int_{\Omega} \left\{\dfrac{1}{m} f(x,u)u - F(x,u)\right\} dx.
$$

\nd Furthermore, by Proposition \ref{lll},
$$
u \in W^{1,\Phi}_{0}(\Omega) \mapsto \int_{\Omega} (\Phi(|\nabla u|)- \frac{1}{m} \phi(|\nabla u|) |\nabla  u|^{2} )dx~~\mbox{  is  w.l.s.c.}
$$

\nd Now using once more  that
$$
u \mapsto \int_{\Omega} F(x, u)dx,~~  u \mapsto \int_{\Omega} f(x, u) u dx~\mbox{are  weakly continuous}
$$we conclude that
\begin{eqnarray*}
c_{\mathcal{N}} &<& \lim_{n \rightarrow \infty} \displaystyle \int_{\Omega} \Phi(|\nabla u_{n}|)- \frac{1}{m} \phi(|\nabla u_{n}|) |\nabla u_{n}|^{2}dx + \displaystyle \int_{\Omega} \left\{\dfrac{1}{m} f(x,u_{n})u_{n} - F(x,u_{n})\right\} \nonumber \\
&=&  \lim_{n \rightarrow \infty} \left\{J(u_{n}) - \dfrac{1}{m} J^{\prime} (u_{n})u_{n} \right\} = c_{\mathcal{N}}, \nonumber \\
\end{eqnarray*}
\nd impossible. Thus  $t = 1$ and $u \in \mathcal{N}$. This finishes the proof. $\hfill{\rule{2mm}{2mm}}$
\vskip.3cm

\nd {\bf Proof of Theorem \ref{th1} (conclusion).}~ Let $(u_{n}) \subset \mathcal{N}$ be a minimizing sequence for $J$ over $\mathcal{N}$. By the proof of  Lemma \ref{pr1} there is
$u \in \mathcal{N} \subset \w$ such that $u_n \rightharpoonup u$.
\vskip.2cm

\nd {\bf Claim.}~~~  $u_n \rightarrow u~\mbox{in}~\w$.
\vskip.2cm

\nd Assume the {\bf Claim} has been proved. Since $J \in C^{1}$ it follows that $J^{\prime}(u_n) \rightarrow J^{\prime}(u)$.
\vskip.1cm

\nd By Lemma \ref{pr1}, $u \in \mathcal{N}$ and
$$
c_{\mathcal{N}} = J(u) = \min_{\mathcal{N}} J > 0.
$$
\nd By Proposition \ref{c11}  the set $\mathcal{N}$ is a ${C}^1$-submanifold of $W^{1,\Phi}_{0}(\Omega)$ so that u is a critical point of $J_{\mid \mathcal{N}}$ and yet by Proposition \ref{c11}  $u$ is a critical point of $J$.
\vskip.2cm

\nd {\bf Proof of the Claim}~ Let $(u_n)$ be the minimizing sequence for $c_{\mathcal{N}}$. Applying the compactness of the  embedding  $W_{0}^{1,\Phi}(\Omega) \hookrightarrow L_\Phi(\Omega)$, it follows  that
$u_{n} \rightarrow u$ in $L_{\Phi}(\Omega)$.
\vskip.1cm

\nd Arguing by contradiction,  there is $\delta > 0$ such that
\begin{equation}\label{r1}
\liminf_{n \rightarrow \infty} \int_{\Omega} \Phi(|\nabla u_{n} - \nabla u|) dx \geq \delta > 0.
\end{equation}
\nd By the  Br\'ezis-Lieb Lemma  for convex functions,  we have
\begin{equation}\label{r2}
\lim_{n \rightarrow \infty} \int_{\Omega} \Phi(|\nabla u_{n}|)  - \Phi(|\nabla u_{n} - \nabla u|) dx = \int_{\Omega} \Phi(|\nabla u|) dx.
\end{equation}
\nd Using \eqref{r1} and \eqref{r2} we infer that
\begin{equation}\label{r3}
 \int_{\Omega} \Phi(|\nabla u|)
 \leq \lim_{n \rightarrow \infty} \int_{\Omega} \Phi(|\nabla u_{n}|) dx - \delta < \lim_{n \rightarrow \infty} \int_{\Omega} \Phi(|\nabla u_{n}|) dx. \nonumber \\
\end{equation}
\nd As a consequence, we get by using the Lebesgue Theorem that
\begin{equation*}
c_{\mathcal{N}} = \lim_{n \rightarrow \infty} J(u_{n}) = \lim_{n \rightarrow \infty} \left\{\int_{\Omega} \Phi(|\nabla u_{n}|) dx - \int_{\Omega} F(x, u_{n}) dx \right\} > J(u)
\end{equation*}
\nd which impossible  because of  $\displaystyle c_{\mathcal{N}} = \lim_{n \rightarrow \infty} J(u_{n}) = J(u)$. Therefore  $u_{n} \rightarrow u$ in $W^{1,\Phi}_{0}(\Omega)$.  This finishes the proof. $\hfill{\rule{2mm}{2mm}}$

\section{\bf \Large Proof of Theorem \ref{th2}}

\nd Consider the two auxiliary  functions
\begin{equation*}
f^{+}(x,t) = \left\{
\begin{array}{l}
f(x,t) \, \,\mbox{if} \,\,t \geq 0, \ \ \\
 0 \,\,\mbox{if} \,\, t < 0 \ \ \\
\end{array}
 \right.
\end{equation*}
and
\begin{equation*}
f^{-}(x,t) = \left\{
\begin{array}{l}
f(x,t) \, \,\mbox{if} \,\,t \leq 0, \ \ \\
 0 \,\,\mbox{if} \,\, t > 0. \ \ \\
\end{array}
 \right.
\end{equation*}
The associated functionals are $J_{\pm} : W^{1,\Phi}_{0}(\Omega)\rightarrow \mathbb{R}$ given by
\begin{equation*}
J_{\pm}(u) = \int_{\Omega} \Phi(|\nabla u|) dx - \int_{\Omega} F^{\pm}(x,u)dx, u \in W^{1,\Phi}_{0}(\Omega)
\end{equation*}
where $F^{\pm}(x,t) = \int_{0}^{t} f^{\pm}(x,s) dx, x \in \Omega, t \in \mathbb{R}$. Recalling definitions in Section 3 we see that  the Nehari manifolds for the functions $f^{+}$ and $f^{-}$ are respectively,
\begin{equation*}
\mathcal{N}^{\pm} = \{ u \in W^{1,\Phi}_{0}(\Omega) \setminus \{0\}~|~  \langle J_{\pm}^{\prime}(u), u\rangle = 0 \}
\end{equation*}
\nd By Proposition \ref{l22}, $\mathcal{N}^{\pm}$ are $C^{1}$-manifolds which are  away from zero. In addition,
\begin{equation*}
c^{\pm} = \inf_{v \in \mathcal{N}^{\pm}} J_{\pm}(v)
\end{equation*}\label{S}
\nd are critical values of $J_{\pm}$. So  we obtain two critical points say $u_{1}, u_{2} \in W^{1,\Phi}_{0}(\Omega) \setminus \{0\} $ such that
$$
J_{+}(u_{1}) = c^{+} > 0~ \mbox{and}~  J_{-}(u_{2}) = c^{-} >0.
$$

\nd Given $u \in \w(\Omega)$ set  $u^{+} = \max\{u,0\},~ u^{-} = \min\{u,0\}$ so that $u = u^{+} + u^{-}$.
\vskip.1cm

\nd  Using $u_{1}^{-}$ as a test function we have
\begin{equation}\label{des-1}
0 \leq \int_{\Omega} \Phi(|\nabla u_{1}^{-}|) \leq \frac{1}{\ell} \int_{\Omega} \phi(|\nabla u_{1}^{-}|)|\nabla u^{-}_{1}|^{2} dx = \frac{1}{\ell} \int_{\Omega} f^{+}(x,u_{1}) u_{1}^{-} dx = 0.
\end{equation}
As a consequence of (\ref{des-1})   one has  $u_{1}^{-} \equiv 0$ and so $u_{1} = u_{1}^{+} \geq 0$ in $\Omega$. Similarly, we also obtain $u_{2} = u_{2}^{-} \leq 0$. Now, since $\mathcal{N}^{\pm}$ are away from zero, it follows that $u_{1}, u_{2} \neq 0$. Hence  by the Maximum Principle,  (cf.  Pucci \& Serrin \cite{pucci}), $u_{1} > 0$ and $u_{2} < 0$ in $\Omega$. For further comments in this regard we  refer  the reader to Carvalho et al \cite{JVMLED}.
\vskip.1cm

\nd We add that the functions $u_{1}, u_{2}$ are critical points of  $J$.  So  $u_{1}, u_{2}$ are also weak solutions of problem \eqref{pi}.
\vskip.1cm

\nd In what follows we shall prove that  problem \eqref{pi} admits at leat one sign changing solution. At first, we define
the nodal Nehari manifold by
\begin{equation*}
\mathcal{N}_{nod} = \{ u \in W^{1,\Phi}_{0}(\Omega) \backslash \{0\}~|~  u^{\pm} \neq 0,  \langle J^{\prime}(u), u^{\pm} \rangle = 0 \}.
\end{equation*}
Consider the nodal level given by
\begin{equation*}
c_{nod} = \displaystyle\inf_{v \in \mathcal{N}_{nod} } J(v).
\end{equation*}
It is easy to verify that any sign changing solution for the problem \eqref{pi} should belong to $\mathcal{N}_{nod}$. Hence it is natural to consider the nodal Nehari manifold in order to ensure the existence of sign changing solutions.
\vskip.1cm

\nd So,  our aim is to  prove that the Nehari manifold $\mathcal{N}_{nod}$ is not empty.
\vskip.1cm

\nd In this regard,   consider the function
$\theta : \mathbb{R}^{+} \times \mathbb{R}^{+} \rightarrow \mathbb{R}$  defined by
\begin{equation*}
\theta(t,s) := J(t u^{+} + s u^{-})
\end{equation*}
where $u \in W^{1,\Phi}_{0}(\Omega)$ and  $u^{\pm} \neq 0$. Given $t, s > 0$ we emphasize  that $\nabla \theta(t ,s) = 0$ if only if $t u^{+} + s u^{-} \in \mathcal{N}_{nod}$. In other words, the critical points of  the function $\theta$ provide us with elements on the Nehari nodal set $\mathcal{N}_{nod}$.
\vskip.1cm

\nd  Given $u \in W^{1,\Phi}_{0}(\Omega)$ with $u^{\pm} \neq 0$,  using once more  the fibering maps, we will prove that
$t u^{+} + s u^{-} \in \mathcal{N}_{nod}$ for some suitable numbers $t , s \in (0, \infty)$. We need a  technical result which we state and prove below.
\vskip.1cm

\begin{prop}\label{est}
Assume  $(\phi_{1}) - (\phi_{3})$ , $(f_{0}) - (f_{3})$. Let $u \in \w $ with  $u^{\pm} \neq 0$.

\nd Then
\begin{equation*}
\lim_{t \rightarrow \infty, s \rightarrow \infty } \dfrac{\theta(t,s)}{t^{m} + s^{m}} = - \infty.
\end{equation*}
\nd and
\begin{equation*}
\lim_{t \rightarrow 0, s \rightarrow 0 } \dfrac{\theta(t,s)}{t^{m} + s^{m}} > 0.
\end{equation*}
\end{prop}
\nd {\bf Proof of Proposition  \ref{est}.} Let $t, s >0$. We have
\begin{eqnarray*}
\theta(t,s)  &=& \int_{\Omega} \Phi(|\nabla t u^{+}|) dx - \int_{\Omega} F(x, t u^{+}) dx \nonumber \\&+& \int_{\Omega} \Phi(|\nabla s u^{-}|) dx - \int_{\Omega} F(x, s u^{-}) dx
\nonumber\\
&\geq& \left(1 - \dfrac{\lambda_1 - \varepsilon}{\lambda_{1}}\right) \int_{\Omega} \Phi(|\nabla t u^{+}|) dx +  C_{\varepsilon}\int_{\Omega} \Psi(t u^{+})dx \nonumber \\
&+& \left(1 - \dfrac{\lambda_1 -  \varepsilon}{\lambda_{1}}\right) \int_{\Omega} \Phi(|\nabla s u^{-}|) dx + C_{\varepsilon}\int_{\Omega} \Psi(s u^{-}) dx. \nonumber\\
\end{eqnarray*}
\nd Applying  Proposition \ref{lema_naru} twice we get to
\begin{eqnarray*}
\theta(t,s) \geq  \left(1 - \dfrac{\lambda_1 - \varepsilon}{\lambda_{1}}\right) t^{m} \int_{\Omega} \Phi(|\nabla u^{+}|) dx - C_{\varepsilon}t^{\ell_{\Psi}}\int_{\Omega} \Psi(u^{+})dx+ \nonumber\\
\left(1 - \dfrac{\lambda_1 - \varepsilon}{\lambda_{1}}\right) s^{m} \int_{\Omega} \Phi(|\nabla u^{-}|) dx - C_{\varepsilon}s^{\ell_{\Psi}}\int_{\Omega} \Psi(u^{-})dx \nonumber\\
\end{eqnarray*}
\nd It follows that
\begin{eqnarray*}
\dfrac{\theta(t,s)}{t^{m} + s^{m}} &\geq&  \left(1 - \dfrac{\lambda_1 - \varepsilon}{\lambda_{1}}\right) \min\left[\int_{\Omega} \Phi(|\nabla u^{+}|) dx , \int_{\Omega} \Phi(|\nabla u^{-}|) dx \right] \nonumber \\
&-& C_{\varepsilon} \dfrac{\max \left(t^{\ell_{\Psi}}, s^{\ell_{\Psi}}\right)}{t^{m} + s^{m}}\int_{\Omega} \Psi(u)dx. \nonumber \\
\end{eqnarray*}
\nd As a consequence,
\begin{equation*}
\lim_{t \rightarrow 0, s \rightarrow 0} \dfrac{\theta(t,s)}{t^{m} + s^{m}} \geq  \left(1 - \dfrac{\lambda_1 - \varepsilon}{\lambda_{1}}\right) \min\left[\int_{\Omega} \Phi(|\nabla u^{+}|) dx , \int_{\Omega} \Phi(|\nabla u^{-}|) dx \right]> 0
\end{equation*}
\nd  for each $\varepsilon > 0$ small enough.
\vskip.2cm

\nd  Next we shall  compute the limit  at infinity of $\theta(t,s)/(t^{m} + s^{m})$. By Proposition \ref{lema_naru},
\begin{equation*}\label{aa1}
\theta(t,s) \leq  t^{m} \int_{\Omega} \Phi(|\nabla u^{+}|) dx + s^{m} \int_{\Omega} \Phi(|\nabla u^{-}|) dx - \int_{\Omega} F(x, t u^{+} + s u^{-}) dx \\
\end{equation*}
\nd which yields
\begin{equation*}
\dfrac{\theta(t,s)}{t^{m} + s^{m}} \leq \max \left(\int_{\Omega} \Phi(|\nabla u^{+}|) dx , \int_{\Omega} \Phi(|\nabla u^{-}|) dx \right)
- \displaystyle \dfrac{1}{t^{m} + s^{m}} \int_{\Omega} F(x, t u^{+} + s u^{-}) dx
\end{equation*}
\nd for $t , s \geq 1$. As a consequence,  we have
\begin{eqnarray}\label{g1}
\displaystyle \lim_{t \rightarrow \infty, s \rightarrow \infty} \dfrac{\theta(t,s)}{t^{m} + s^{m}} &\leq&  \max\left(\int_{\Omega} \Phi(|\nabla u^{+}|) dx , \int_{\Omega} \Phi(|\nabla u^{-}|) dx \right) \nonumber \\ &-& \liminf_{t \rightarrow \infty, s \rightarrow \infty} \int_{\Omega} \dfrac{F(x, t u^{+} + s u^{-})}{t^{m} + s^{m}} dx. \nonumber \\
\end{eqnarray}
\nd We have,
\begin{equation*}
\liminf_{t \rightarrow \infty, s \rightarrow \infty} \int_{\Omega} \dfrac{F(x, t u^{+} + s u^{-})}{t^{m} + s^{m}} dx \geq
\liminf_{t \rightarrow \infty, s \rightarrow \infty} \int_{\Omega} \dfrac{F(x, t u^{+} + s u^{-})}{|t u^{+} + s u^{-}|^{m}} \dfrac{|t u^{+} + s u^{-}|^{m}} {t^{m} + s^{m}} dx
\end{equation*}
\nd  By  Fatou's Lemma and $(f_{3})$  it follows that\\
\begin{eqnarray}\label{g2}
  \infty & = & \liminf_{t \rightarrow \infty, s \rightarrow \infty} \int_{\Omega} \dfrac{F(x, t u^{+} + s u^{-})}{|t u^{+} + s u^{-}|^{m}} \min(|u^{+}|^{m},|u^{-}|^{m}) dx\nonumber \\
&\leq&\liminf_{t \rightarrow \infty, s \rightarrow \infty} \int_{\Omega} \dfrac{F(x, t u^{+} + s u^{-})}{t^{m} + s^{m}} dx.
\end{eqnarray}
\nd It follows from \eqref{g1} and \eqref{g2} that
\begin{equation*}
\lim_{t \rightarrow \infty, s \rightarrow \infty} \dfrac{\theta(t,s)}{t^{m} + s^{m}} = -\infty.
\end{equation*}
This ends the proof.     $\hfill{\rule{2mm}{2mm}}$
\vskip.2cm

\begin{prop}\label{Nehari not empty}
Assume $(\phi_{1}) - (\phi_{3})$, $(f_{0}) - (f_{3})$. Let $u \in W^{1,\Phi}_{0}(\Omega)$ such that   $u^{\pm} \neq 0$. Then  there exist uniquely determined $t, s \in (0, \infty)$ such that $t u^{+} + s u^{-} \in \mathcal{N}_{nod}$.
\end{prop}

\nd {\bf Proof of Proposition  \ref{Nehari not empty}.}  By Proposition \ref{est} we infer that
\begin{equation*}
\lim_{t \rightarrow \infty, s \rightarrow \infty } \dfrac{\theta(t,s)}{t^{m} + s^{m}} = - \infty.
\end{equation*}
and
\begin{equation*}
\lim_{t \rightarrow 0, s \rightarrow 0 } \dfrac{\theta(t,s)}{t^{m} + s^{m}} > 0.
\end{equation*}
\nd Since  $\theta(t,s)$ is continuous,  there exist $t_0,s_0 \in (0,\infty)$ such that
$$
\displaystyle\max_{t \geq 0, s \geq 0 } \theta(t,s) = \theta(t_{0},s_{0}).
$$
By  Proposition \ref{l1} and
$$
\frac{\partial \theta(t_0,s_0)}{\partial t}=\gamma^\prime_{u^+}(t_0)=0,~~~
\frac{\partial \theta(t_0,s_0)}{\partial s}=\gamma^\prime_{u^+}(s_0)=0
$$
\nd it follows that $t_0, s_0$ are uniquely determined. By an easy argument  one shows that  $t_{0} u^{+} + s_{0} u^{-} \in \mathcal{N}_{nod}$. This ends the proof.  $\hfill{\rule{2mm}{2mm}}$
\vskip.2cm

\begin{prop}\label{Nehari not empty-1}
Assume  $(\phi_{1}) - (\phi_{3})$, $(f_{0}) - (f_{3})$. Let  $(u_{n}) \subseteq \mathcal{N}_{nod}$    be a minimizing sequence for $J$. Then $(u_n)$ is bounded in $\w$.
\end{prop}
\nd {\bf Proof of Proposition  \ref{Nehari not empty-1}.} The proof  follows the same lines as in the proof   of Lemma \ref{l33}.  Since
$u_{n}  = u_{n}^{+} + u_{n}^{-}$  it suffices to show that  $(u_{n}^{+})$ and $(u_{n}^{-})$  are bounded in $\w$. We skip the details. $\hfill{\rule{2mm}{2mm}}$

\begin{prop}\label{n1}
Assume $(\phi_{1}) - (\phi_{3})$, $(f_{0}) - (f_{3})$. Then there is $u \in \mathcal{N}_{nod}$  such  that $J(u) = c_{nod}$.
\end{prop}
\vskip.2cm

\nd  {\bf Proof of Proposition  \ref{n1}.}   Let $(u_{n})  \subseteq \mathcal{N}_{nod}$ be a minimizing sequence for $J$ over $ \mathcal{N}_{nod}$, that is
$$
J(u_n) \rightarrow  c_{nod}.
$$
\nd We split  $u_n$ as  $u_{n}  = u_{n}^{+} + u_{n}^{-}$. By  Proposition  \ref{Nehari not empty-1},  $(u_n)$ is bounded in the reflexive space $\w(\Omega)$ and as a consequence,  both  $(u^+_{n})$ and  $(u^-_{n})$ are also bounded in  $\w(\Omega)$ and by the very definition of   $ \mathcal{N}_{nod}$,  $u_{n}^{\pm} \neq 0$. Thus,
$$
u_n \rightharpoonup u,~ u_n^{+}  \rightharpoonup u^{+},~ u_n^{-}  \rightharpoonup u^{-}~\mbox{in}~\w(\Omega),
$$
$$
u^{+} \cdot u^{-} = 0,~\mbox{supp}(u_n^{+}) \cap \mbox{supp}(u_n^{-}) = \emptyset,
$$
$$
\langle J^{\prime}( u_n),  u_n \rangle = \langle J^{\prime}( u_n^{+}),  u_n^{+} \rangle = \langle J^{\prime}( u_n^{-}),  u_n^{-} \rangle = 0
$$
\nd so that, in particular,  $(u_{n}^{\pm})  \subseteq   \mathcal{N} $.

\vskip.2cm

\nd {\bf Claim.}~~~~  $u^{+},~  u^{-} \neq  0$.
\vskip.1cm

\nd Indeed, assume on the contrary that  $u^{+} = 0$. Using the fact that  $u_{n}^{+}  \in \mathcal{N}$ we have
$$
\begin{array}{lll}
\displaystyle  0 \leq \displaystyle \int_{\Omega} \Phi(|\nabla u_{n}^{+} |)dx &\leq& \displaystyle \frac{1}{\ell} \int_{\Omega} \phi(|\nabla u_{n}^{+}|)|\nabla u^{+}_{n}|^{2} dx   \\
\\
&= &  \displaystyle  \frac{1}{\ell} \int_{\Omega} f(x,u_{n}^{+}) u_{n}^{+} dx \\
\\
&=&  \frac1 \ell~ o_n(1).
\end{array}
$$
\nd It follows that
$$
\displaystyle \min \{ \Vert u_n^{+} \Vert^{\ell},~\Vert u_n^{+} \Vert^{m} \} \leq \int_{\Omega} \Phi(|\nabla u_{n}^{+} |) dx = o_n(1)
$$
\nd and so $ \Vert u_n^{+} \Vert \rightarrow 0.$  Assuming that $u^{-} = 0$ we find  in a similar way  that   $ \Vert u_n^{-} \Vert \rightarrow 0$. By Proposition \ref{l22}  and the fact that $u_n^{\pm}  \rightharpoonup u^{\pm}~\mbox{in}~\w(\Omega)$  we have
$$
\displaystyle 0 < C  \leq \Vert u^{\pm} \Vert \leq  \liminf \Vert u_n^{\pm} \Vert
$$
\nd which is impossible. This proves the {\bf Claim.}
\vskip.1cm

\nd By Proposition \ref{Nehari not empty}  there exist $t, s \in (0,\infty)$ such that $t u^{+} + s u^{-} \in \mathcal{N}_{nod}$. Set
$$
\widetilde{u} = t u^{+} + s u^{-}.
$$
\nd At this point we have
\begin{equation}\label{jv-1}
t u^{+}, s u^{-} \neq 0,~~ t u^{+} \cdot s u^{-} = 0,~~ \mbox{supp}(t u^{+}) \cap \mbox{supp}(su^{-}) = \emptyset.
\end{equation}
\vskip.1cm

\nd Using the fact that  $\widetilde{u} = t u^{+} + s u^{-} \in \mathcal{N}_{nod}$ and (\ref{jv-1}) we have
$$
0 = \langle J^{\prime}(\widetilde{u}), {\widetilde{u}}^{\pm} \rangle = \langle J^{\prime}(t {{u}}^{+}), t {{u}}^{+} \rangle = \langle J^{\prime}(s {{u}}^{-}), s {{u}}^{-} \rangle.
$$
\nd  On the other hand, using Proposition \ref{lll} (i),
$$
\langle J^{\prime}( u^{+}),  u^{+} \rangle \leq \liminf  \langle J^{\prime}( u_n^{+}),  u_n^{+} \rangle = 0
$$
\nd and
$$
\langle J^{\prime}( u^{-}),  u^{-} \rangle \leq \liminf  \langle J^{\prime}( u_n^{-}),  u_n^{-} \rangle = 0.
$$
\nd As a consequence of the inequalities just above, $0 < t, s \leq 1$,  because if otherwise, one of $t,s$, say $t > 1$ we would have $\langle J^{\prime}(\tau u^{+}), \tau u^{+} \rangle >0$ for each $\tau \in (0,t)$. But  since  $\gamma_{u^{+}}^{\prime}(1) \leq 0$  we have a contradiction.
\nd Now,
$$
\begin{array}{lll}
\displaystyle  c_{nod} \leq \displaystyle J(\widetilde{u})  = J( t u^{+} + s u^{-})  &=&  J( t u^{+}) +  J(s u^{-})\\
\\
&\leq &  \displaystyle  \liminf J(tu_n^{+})+ \liminf J(su_n^{-})  \\
\\
&\leq&   \liminf (J(tu_n^{+}+su_n^{-})) \\
\\
&\leq&  \displaystyle \liminf \max_{t,s >0} \theta(t,s)\\
\end{array}
$$
\nd But, by the proof of Proposition \ref{Nehari not empty},  we have $\displaystyle\max_{t>0,s>0}\theta(t,s)=\theta(1,1)$. Hence,
$$
\begin{array}{lll}
\displaystyle  c_{nod} &\leq& \displaystyle \liminf \theta(1,1) \\
\\
&=& \liminf  J(  u_n^{+} +  u_n^{-}) = c_{nod}.
\end{array}
$$
\nd This ends the proof of Proposition  \ref{n1}.  $\hfill{\rule{2mm}{2mm}}$

\vskip.1cm

\nd Next we shall prove that any minimizer for $J$ on $\mathcal{N}_{nod}$ is a free critical point.

\begin{prop}\label{n2}
Assume $(\phi_{1}) - (\phi_{3})$, $(f_{0}) - (f_{3})$.  Let $u \in \mathcal{N}_{nod}$ be a minimizer for $J$ over $\mathcal{N}_{nod}$, that is   $J(u) = c_{nod}$. Then u is a free critical point for $J$.
\end{prop}

\nd {\bf Proof of Proposition  \ref{n2}.}  We shall adapt arguments  by Szulkin \verb"&" Weth  in \cite{SW}  to the case of the $\Phi$-Laplace operator. For the reader's convenience, we shall provide a few  details.

\nd Remind  that $\gamma_{u^+}(1)$ and  $\gamma_{u^-}(1)$ are the unique maximum values for the functions $\gamma_{u^+},\gamma_{u^-}$, respectively. We have
\begin{equation}\label{descnod}
    \begin{array}{lll}
    J(su^++tu^-) & = & J(su^+)+J(tu^-)\\
\\
&<& J(u^+)+J(u^-)\\
    \\
& = & J(u)\\
\\
&=& c_{nod},~~s,t\geq 0,~~s,t\neq1.
    \end{array}
\end{equation}
In order to show that  $u \in  \mathcal{N}_{nod}$ is a critical point for $J$, assume by the way of contradiction that  $J'(u)\neq 0$.
\vskip.1cm

\nd Then there exist $\delta,\mu>0$ such that
$$
\|J'(v)\|\geq \mu~~ \mbox{whenever}~~ \|v-u\|\leq 3\delta,~~ v\in W^{1,\Phi}_0(\Omega).
$$
 \nd Let  $D= D_{0}\times D_{0}$ where $D_{0} = (\frac12,\frac23)$. Consider   the continuous function  $g:\overline{D}\rightarrow W^{1,\Phi}_0(\Omega)$  defined by $g(s,t)=tu^++su^-$. Note that $t \rightarrow \gamma_{u^+}(t)$ and $t \rightarrow \gamma_{u^-}(t)$ are strictly increasing functions on $[0,1]$. It follows that
$$
J(g(s,t))=J(su^++tu^-)=c_{nod}~ \mbox{if and only if}~ t=s=1.
$$
Moreover, due the estimates  in (\ref{descnod}), we have  $J(g(s,t))< c_{nod}$ for $t , s > 1$. Hence
$$
\beta:=\max_{\partial D}J(g(s,t))< c_{nod}.
$$
\nd By the Deformation Lemma (cf. \cite[Lemma 2.3]{Willem}) with $\epsilon := \min \left\{ \frac{c_{nod}-\beta}{4}, \frac{\mu\delta}{8}\right\}$ there is $\eta\in C\left([0,1]\times W^{1,\Phi}_0(\Omega),W^{1,\Phi}_0(\Omega)\right)$  such  that
\begin{enumerate}
  \item $\eta(1,v)=v$ if $v\not\in J^{-1}([c_{nod}-2\epsilon,c_{nod}+2\epsilon])$,

 \item $J(\eta(1,v))\leq c_{nod}-\epsilon,$ for each $v\in\w$\\
 such that $\|v-u\|\leq \delta$ and $J(v)\leq c_{nod}+\epsilon$,

\item $J(\eta(1,v))\leq J(v)$ for each $v\in \w(\Omega)$.
\end{enumerate}
\nd Since  the maximum value of $J\circ g$ is achieved at $(s,t)=(1,1)$ it follows by  item (3) just above that
$$
\begin{array}{lll}
\displaystyle\max_{(s,t)\in D}J(\eta(1,g(s,t))) & \leq & \displaystyle\max_{(s,t)\in D}J(g(s,t))\\
\\
& < & \displaystyle c_{nod}\\
\\
\displaystyle &=& \displaystyle \max_{\{(s,t)\in [0,1]\times[0,1]\}} J(g(s,t)),
\end{array}
$$
\nd so that   $\eta(1,g(s,t))\not\in \mathcal{N}_{nod}$ for each $(s,t)\in D$.
\vskip.2cm

\nd Now we claim that there is $(s_0,t_0)\in D$ such that $\eta (1,g(s_0,t_0))\in \mathcal{N}_{nod}$ which will  lead to a constradiction, showing that  $u$ is a critical point of $J$.
\vskip.1cm

\nd The proof of the claim will be achieved  using degree theory. Consider the maps
$$
h(s,t):=\eta(1,g(s,t))=\eta(1, su^++tu^-),
$$
$$\Psi_0(s,t)=\left(\langle J'(su^+),u^+\rangle,\langle J'(su^-),u^-\rangle\right)$$
and
$$
\Psi_1(s,t):=\left(\frac1s\langle J'(h^+(s,t)),h^+(s,t)\rangle,\frac1t\langle J'(h^-(s,t)),h^-(s,t)\rangle\right).
$$
\nd We claim that $\deg(\Psi_0, D, (0,0))=1$. Indeed,
$$
\langle J'(su^\pm), u^\pm \rangle >0,~ s\in(0,1),~\mbox{due the fact that}~ \gamma_{u^\pm}(s)>0~\mbox{for}~s\in (0,1)
$$
\nd  and
$$\langle J'(su^\pm),u^\pm\rangle<0,~ s\in(1,+\infty)~~\mbox{because}~~ \gamma_{u^\pm}(s)<0~\mbox{for}~s \in (1,+\infty).$$
\nd  Now, consider the homotopy
$$
H(\tau, s)=(1-\tau)L(s)+\tau s,
$$
\nd where $L(s):=\langle J'(su^+),u^+\rangle.$
\vskip.1cm

\nd Since $H(\tau,s)\neq 0$ for each $\tau\in [0,1]$ and $s\in \{\frac12,\frac23\}$, we have by the  homotopy invariance property that
$$\deg(L,D_{0},0)=\deg(\mbox{Id}_{\mathbb{R}},D_{0} ,0)=1.$$
By product formula for the degree  the claim holds.
\vskip.1cm

\nd Note that
$$
J(g(s,t))\leq \beta =\max_{\partial D}J\circ g\leq c_{nod}-4\epsilon<c_{nod}-2\epsilon.
$$
Hence $g(s,t)\not\in J^{-1}([c_{nod}-2\epsilon,c_{nod}+2\epsilon])$. So,  by item (1) above we see that
$$h(s,t):=\eta(1,g(s,t))=g(s,t),~(s,t)\in\partial D.$$
\nd   Now,  note that
$$
\begin{array}{lll}
\Psi_1(s,t) & = & \displaystyle\left(\frac1s \langle J'(g^+(s,t)),g^+(s,t)\rangle,\langle\frac1t J'(g^-(s,t)),g^-(s,t)\rangle\right)\\
           \\
 & = & \displaystyle\left( \langle J'(su^+),u^+\rangle, \langle J'(tu^-),u^-\rangle\right)\\
\\
&=& \Psi_0(s,t),~~(s,t)\in \partial D,
\end{array}
$$
\nd so that   $\Psi_1=\Psi_0$ on $\partial D$. This implies that
$$\deg(\Psi_1, D, (0,0))=\deg(\Psi_0, D, (0,0))=1.$$
As a consequence there is $(s_0,t_0)\in D$ such that $\Psi_1(s_0,t_0)=(0,0)$. This is equivalent to
$$
\langle J'(h^\pm(s_0,t_0)),h^\pm(s_0,t_0)\rangle =0,
$$
\nd showing that $h^\pm(s_0,t_0)\in \mathcal{N}$.  Therefore $\eta(1,g(s_0,t_0))=h(s_0,t_0)\in\mathcal{N}_{nod}$ which is a contradiction.
\vskip.1cm

\nd  Applying  Propositions \ref{n1} and \ref{n2}  ends the proof of Theorem \ref{th2}.
  $\hfill{\rule{2mm}{2mm}}$

\begin{flushright}
{\small F. J. S. A. Corr\^ea}\\
{\small  email:   fjsacorrea@gmail.com }
\\

{\small  address:}
\\
\scriptsize{ Universidade Federal de Campina Grande\\
 Unidade Acad\^emica de Matem\'atica\\
58.429-900 - Campina Grande - PB - Brazil}\\
\vskip.3cm

{\small  M. L. M. Carvalho}\\
{\small  email: marcos\_leandro\_carvalho@ufg.br},\\
\vskip.3cm

{\small  Jose V. A. Goncalves}\\
{\small  email: goncalvesjva@ufg.br}\\
\nd and\\
\vskip.2cm

{\small    E. D. Silva}\\
{\small email: edcarlos@ufg.br}
\vskip.2cm

{\bf  address:}
\vskip.2cm

\scriptsize{  Universidade Federal de Goi\'as\\
   Instituto de Matem\'atica e Estat\'istica\\
   74001-970 Goi\^ania, GO - Brasil}
\end{flushright}

\end{document}